\documentclass[a4paper,10pt,reqno]{amsart}
\usepackage[a4paper,tmargin=40mm,bmargin=35mm,hmargin=30mm]{geometry}
\usepackage[T1]{fontenc}
\usepackage{lmodern}
\usepackage{amsmath}
\usepackage{amssymb}
\usepackage{amsthm}
\usepackage{stmaryrd}
\usepackage{tikz}
\usepackage[all]{xy}
\usepackage{booktabs}
\usepackage{multirow}
\usepackage{multicol}
\usepackage{here}
\usepackage{aliascnt}
\usepackage{hyperref}
\usepackage[noabbrev]{cleveref}

\newcommand{\NewTheorem}[2]{
	\newaliascnt{#1}{TheoremEnvironment}
	\newtheorem{#1}[#1]{#1}
	\aliascntresetthe{#1}
	\crefname{#1}{#1}{#2}
	\Crefname{#1}{#1}{#2}
}

\theoremstyle{definition}

\NewTheorem{Definition}{Definitions}
\NewTheorem{Remark}{Remarks}
\NewTheorem{Axiom}{Axioms}
\NewTheorem{Example}{Examples}
\NewTheorem{Observation}{Observations}
\NewTheorem{Convention}{Conventions}
\NewTheorem{Notation}{Notations}
\NewTheorem{Setting}{Settings}
\NewTheorem{Question}{Questions}
\NewTheorem{Answer}{Answers}
\NewTheorem{Conjecture}{Conjectures}
\NewTheorem{Problem}{Problems}
\NewTheorem{Solution}{Solutions}
\NewTheorem{Goal}{Goals}
\NewTheorem{Comment}{Comments}
\NewTheorem{Aim}{Aims}
\NewTheorem{Caution}{Cautions}
\NewTheorem{Exercise}{Exercises}
\NewTheorem{Examples}{Examples}
\theoremstyle{plain}
\NewTheorem{Proposition}{Propositions}
\NewTheorem{Lemma}{Lemmas}
\NewTheorem{Theorem}{Theorems}
\NewTheorem{Corollary}{Corollaries}

\crefname{enumi}{}{}
\Crefname{enumi}{}{}
\creflabelformat{enumi}{(#2#1#3)}
\crefname{enumii}{}{}
\Crefname{enumii}{}{}
\creflabelformat{enumii}{(#2#1#3)}
\crefname{enumiii}{}{}
\Crefname{enumiii}{}{}
\creflabelformat{enumiii}{(#2#1#3)}

\makeatletter
\renewcommand{\p@enumii}{}
\renewcommand{\p@enumiii}{}
\makeatother

\numberwithin{equation}{section}
\crefname{equation}{}{}
\Crefname{equation}{}{}
\creflabelformat{equation}{(#2#1#3)}

\newcommand{\SwapSymbols}[1]{
	\expandafter\let\expandafter\temporarysymbol\csname #1\endcsname
	\expandafter\let\csname #1\expandafter\endcsname\csname var#1\endcsname
	\expandafter\let\csname var#1\endcsname\temporarysymbol
}

\SwapSymbols{epsilon}
\SwapSymbols{phi}
\SwapSymbols{Gamma}
\SwapSymbols{Delta}
\SwapSymbols{Theta}
\SwapSymbols{Lambda}
\SwapSymbols{Xi}
\SwapSymbols{Pi}
\SwapSymbols{Sigma}
\SwapSymbols{Upsilon}
\SwapSymbols{Phi}
\SwapSymbols{Psi}
\SwapSymbols{Omega}

\newcommand{\cA}{\mathcal{A}}
\newcommand{\cB}{\mathcal{B}}

\newcommand{\cE}{\mathcal{E}}
\newcommand{\cF}{\mathcal{F}}
\newcommand{\cG}{\mathcal{G}}

\newcommand{\cI}{\mathcal{I}}
\newcommand{\cJ}{\mathcal{J}}

\newcommand{\cX}{\mathcal{X}}

\newcommand{\To}{\longrightarrow}

\DeclareMathOperator{\Hom}{Hom}

\DeclareMathOperator{\Tor}{Tor}
\DeclareMathOperator{\id}{id}
\DeclareMathOperator{\Gid}{Gid}
\DeclareMathOperator{\RHom}{{\bf R}Hom}

\DeclareMathOperator{\Gfd}{Gfd}

\DeclareMathOperator{\Ker}{Ker}
\DeclareMathOperator{\Coker}{Coker}
\let\Im\relax
\DeclareMathOperator{\Im}{Im}

\DeclareMathOperator{\Spec}{Spec}

\DeclareMathOperator{\Ass}{Ass}

\DeclareMathOperator{\height}{ht}

\title{Gorenstein injective filtrations over rings with dualizing complexes}

\author{Reza Sazeedeh}

\address{Department of Mathematics, Urmia University, P.O.Box: 165, Urmia, Iran}
\email{rsazeedeh@ipm.ir and r.sazeedeh@urmia.ac.ir}

\begin{document}

\begin{abstract}
Let $R$ be a commutative noetherian ring. Enochs and Huang [EH] proved that over a Gorenstein ring of Krull dimension $d$, every Gorenstein injective module admits a finite filtration of Gorenstein injective submodules. In this paper, we extend this result to rings admitting a dualizing complex and we provide such filtrations using Auslander categories and  section functors. 
\end{abstract}

\maketitle
\tableofcontents

\section{Introduction}
Throughout this paper, $R$ is a commutative noetherian ring with identity admitting a dualizing complex $D$. We recall from [EJ] that an $R$-module $G$ is Gorenstein injective if there exists an exact sequence of injective $R$-modules 
$$\cI:\dots\To E_1\To E_0\To E_{-1}\To \dots$$ such that  $\Hom_R(E,\cI)$ is exact for all injective $R$-modules $E$ and $G=\Ker(E_0\To E_{-1})$. 

 Matlis [M] proved that over a noetherian ring, every injective module is the direct sum of  indecomposable injective modules. A result analogous to Matlis's theorem for Gorenstein injective modules is the existence of finite filtrations by Gorenstein injective submodules. This issue was initially pursued by Enochs and Huang [EH] over Gorenstein rings of finite Krull dimension and they asked whether such filtrations exist when $R$ is a Cohen-Macaulay ring with a dualizing module. Feickert and Sather-Wagstaff [FS] provided a positive response and subsequently raised the same question when $R$ admits a dualizing complex. 
 
 In this paper, in two different ways, we show that if $R$ has a dualizing complex $D$, then such filtrations exist. By setting $\sup D=0$ and $X_k=\{\frak p\in \Spec R|\hspace{0.1cm}{\rm ht}\frak p-\sup D_{\frak p}=k\}$ for every integer $k$, the first method using Auslander categories enables us to prove the following theorem. 

\begin{Theorem}\label{1.1}
Let $G$ be a Gorenstein injective $R$-module. Then $G$ has a finite filtration of Gorenstein injective submodules $$0=G_{d+1}\subset G_{d}\subset\dots \subset G_{1}\subset G_{0}=G$$ such that $G_{k}/G_{k+1}\cong \bigoplus_{\frak p\in X_{k}}\Tor_{k}^R(E(R/\frak p),\RHom_R(D,G))$ is Gorenstein injective for each $0\leq k\leq d$. Furthermore, such filtrations and direct sum decompositions are unique and functorial in $G$. 
\end{Theorem}

When $R$ is a Gorenstein ring, substituting the  dualizing complex $D$ with $R$ recovers  the filtration in \cref{1.1} as identical to that in [EH, Theorem 3.1]. Similarly, for a Cohen-Macaulay ring equipped with a dualizing module $D$ (instead of a dualizing complex $D$), the filtration aligns with the one constructed in [FS, Theorem 4.2]. 

 In a separate approach, when $\dim R=d$, we construct finite filtrations for Gorenstein injective modules using section functors. Unlike \cref{1.1}, where the filtration are characterized via quotients $G_k/G_{k+1}$, here we explicitly build each submodule $G_k$ in the filtration.  This distinction shows the significance of the new filtrtions and their complementary role in analyzing Gorenstein injective structures. The following theorem formulates this filtration:
 
  \begin{Theorem}\label{1.2}
Let $G$ be a Gorenstein injective module and let $Y_k=\{\frak p\in\Spec R|\hspace{0.1cm} \height \frak p=k\}$ for each $k\geq 0$. Then $G$ has a finite filtration  of Gorenstein injective submodules $$0=G_{d+1}\subset G_d\subset G_{d-1}\subset\dots \subset G_1\subset G_0=G$$ such that $G_k= \Sigma_{\frak p\in Y_k}\Gamma_{\frak p}(G)$ and $G_k/G_{k+1}=\Sigma_{\frak p\in Y_k}\Gamma_{\frak p}(G/G_{k+1})$ for each $0\leq k\leq d$. Furthermore, such filtrations and sums  are unique and functorial in $G$. 
\end{Theorem}
Finally, we show that the sums defining $G_k/G_{k+1}$ in \cref{1.2} are direct. This establishes that the filtrations in \cref{1.1,1.2} coincide. To be more precise, we have the following theorem.
\begin{Theorem}
Let $G$ be a Gorenstein injective $R$-module with the the same filtration  as in \cref{1.2}. Then $G_k/G_{k+1}=\bigoplus_{\frak p\in Y_k}\Gamma_{\frak p}(G/G_{k+1})$ for every integer $k$. Furthermore, the filtrations in \cref{1.1,1.2} are the same.
\end{Theorem}
 
\medskip

\section{Filtration for Gorenstein injective modules}
At the beginning of this section, we briefly list some necessary basic notations and definitions. 

For any complex $X=\dots\To X_{i+1}\stackrel{\partial_{i+1}}\To X_i\stackrel{\partial_i}\To X_{i-1}\To\dots $ of $R$-modules,
 the {\it infimum} inf$X$ and the {\it supremum} sup$X$ of $X$ are the (possible infinite) numbers $\inf\{n\in\mathbb{Z}|\hspace{0.1cm} H_n(X)\neq 0\}$ and $\sup\{n\in\mathbb{Z}|\hspace{0.1cm} H_n(X)\neq 0\}.$ The class of  $R$-complexes is denoted by $C(R)$. An $R$-complex $X$ is right (left) {\it bounded} if $X_n=0$  for all $n\ll 0$ ($n\gg 0$) and $X$ is bounded if it is right and left bounded. The category of bounded $R$-complexes is denoted by $C_b(R)$. For each $n\in\mathbb{Z}$, the $n$-{\it th shifted} of $X$ is the complex $\Sigma^nX$ given by $(\Sigma^nX)_i=X_{i-n}$ and $\partial_i^{\Sigma^nX}=(-1)^n\partial_{i-n}$. We use the notation $D(R)$ for the derived category of $R$, and $D_b(R)$ for full subcategory of homologically right and left bounded $R$-complexes. We denote by $D_b^f(R)$, a full subcategory of $D_b(R)$ consisting of complexes $X$ such that $H(X)$ is finitely generated in each degree. The symbol $\simeq$ shows isomorphism in $D(R)$, that is the quasi-isomorphism in $C(R)$.    

\medskip

\begin{Definition} A complex $D\in D_b^f(R)$ is {\it dualizing} for $R$ if it has finite injective dimension and the natural homothety morphism $\cX_D^R:R\To \RHom_R(D,D)$ is an isomorphism in $D(R)$. If $(R,\frak m)$ is a local ring, then the dualizing complex $D$ is {\it normalized} if $\sup D=\dim R.$ 
 \end{Definition}
 
 \medskip
 
 The {\it Auslander categories} $\cA(R)$ and $\cB(R)$ with respect to the dualizing complex $D$ are full subcategories  of $D_b(R)$ defined as:
 
 \medskip
 $\cA(R)=\{X\in D_b(R)|\hspace{0.1cm}\eta_X:X\stackrel{\simeq}\To\RHom_R(D,D\otimes_R^{\bf L}X)$ in $D(R)$ and $D\otimes_R^{\bf L}X\in D_b(R)\}.$

   and

$\cB(R)=\{Y\in D_b(R)|\hspace{0.1cm}\epsilon_Y: D\otimes_R^{\bf L}\RHom(D,Y)\stackrel{\simeq}\To Y$ in $D(R)$ and $\RHom_R(D,Y)\in D_b(R)\}.$

 \medskip
 \begin{Definition} For any complex $M\in D_b(R)$, the {\it Gorenstein injective dimension} of $M$, denoted by $\Gid_RM$ is

 $\Gid_RM=\inf\{\sup\{l\in\mathbb{Z}|\hspace{0.1cm} G_{-l}\neq 0\}|\hspace{0.13cm}G\in C_{b}(R)$ is isomorphic to $M$ in $D(R)$
 
  \hspace{5.5cm} and every $G_l$ is Gorenstein injective $\}$.
  
  The {\it Gorenstein flat dimension} of $M$, denoted by $\Gfd_RM$ is defined dually.
  \end{Definition} 
  \medskip
  The major techniques used in this paper are within the framework of the derived categories and local cohomology. For more study on these subjects, we refer the reader to the text books [BS, CFH, H]. 
  
 It turns out from the proof of [S1, Theorem 3.1] that if $\frak a$ is an ideal of a commutative noetherian ring $R$ (without additional conditions on $R$) and $G$ is a Gorenstein injective $R$-module, then $H_{\frak a}^i(M)=0$ for all $i>0$. We further have the following lemma.

\begin{Lemma}\label{secgor}
Let $\frak a$ be an ideal of $R$ and $G$ be a Gorenstein injective $R$-module. Then $H_{\frak a}^i(G)=0$ for all $i>0$. Furthermore, $\Gamma_{\frak a}(G)$ is Gorenstein injective.
\end{Lemma}
\begin{proof}
 See  [S1, Theorem 3.1] and [S2, Theorem 3.2].
\end{proof}
\medskip

First, we establish several preliminary results  which will be used in the proof of the main theorem. 
\begin{Lemma}\label{zero}
Let $\frak p$ be a prime ideal of $R$ and $G$ be a Gorestein injective $R$-module. Then for all $i\neq \dim R_{\frak p}-\sup D_{\frak p}$, we have $\Tor_i^R(E(R/\frak p),{\bf R}\Hom_R(D,G))=0.$  
\end{Lemma}
\begin{proof}
As $E(R/\frak p)$ is an injective $R_{\frak p}$-module, we have $$\Tor_i^R(E(R/\frak p),{\bf R}\Hom_R(D,G))\cong \Tor_i^{R_{\frak p}}(E(R/\frak p),{\bf R}\Hom_{R_{\frak p}}(D_{\frak p},G_{\frak p})).$$
 We observe that $D_{\frak p}$ is a dualizing complex for $R_{\frak p}$ and by [CFH1, Proposition 5.5], the $R_{\frak p}$-module $G_{\frak p}$ is Gorenstein injective.  Therefore, we may assume that $\frak p=\frak m$ and $(R,\frak m)$ is a local ring with the residue field $k=R/\frak m$. Taking
  $t=\dim R-\sup D$, the shifted complex $\Sigma^{t}D$ is a normalized dualizing complex and so by [CF, Proposition 7.1.11], we have $\id_R \Sigma^{t}D=0$.   
Since $G$ is Gorenstein injective, we have $\RHom_R(D,G)\in\cA(R)$ and so it follows from [CFH1, Lemma 4.5] that
$$\sup(E(R/\frak m)\otimes_R^{\bf L}\RHom_R(D,G))\leq \id_R D+\sup(D\otimes_R^{\bf L}\RHom_R(D,G))$$$$=\id_RD+\sup G=\id_RD=\id_R(\Sigma^{t}D)+t=t.$$ We now prove that $\inf(E(R/\frak m)\otimes_R^{\bf L}\RHom_R(D,G))\geq t$. Since $D$ has finite injective dimension, by [R, Theorem 2.4, p.42], it has a bounded minimal injective resolution $$\cI:0\To I_m\To\dots\To I_{n-1}\To I_n\To 0.$$ It follows from [CF, Proposition 7.1.11] that  rank$_k(H_{-t}(\RHom_R(k,D)))=1$ and for each $i\neq -t$,$$H_i(\RHom_R(k,D))=H_{i+t}(\RHom_R(k,\Sigma^{t}D))=0$$ which means that only $I_{-t}$ contains one copy of $E(R/\frak m)$. On the other hand, for each $\frak q\in\Spec R$ with $\frak q\subsetneq \frak m$ and any $R$-module $X$, taking $r\in\frak m\setminus\frak q$, the linear map $E(R/\frak q)\stackrel{r}\To E(R/\frak q)$ is isomorphism so that $E(R/\frak m)\otimes_R\Hom_R(E(R/\frak q),X)\stackrel{r.}\To E(R/\frak m)\otimes_R\Hom_R(E(R/\frak q),X)$ is isomorphism. Thus the fact that any element of $E(R/\frak m)\otimes_R\Hom_R(E(R/\frak q),X)$ is annihilated by some power of $r$ forces $E(R/\frak m)\otimes_R\Hom_R(E(R/\frak q),X)=0.$ Given an injective resolution $\cE:0\To E_0\To E_{-1}\To \dots $ of $G$, we have $$E(R/\frak m)\otimes_R^{\bf L}\RHom_R(D,G)\simeq E(R/\frak m)\otimes_R\Hom_R(\cI,\cE). $$ Then, by the previous argument, for each $n\in\mathbb{Z}$, we have $$(E(R/\frak m)\otimes_R\Hom_R(\cI,\cE))_n=E(R/\frak m)\otimes_R\Hom_R(I_{-t},E_{n-t})$$$$=E(R/\frak m)\otimes_R\Hom_R(E(R/\frak m),E_{n-t})=(\Sigma^t E(R/\frak m)\otimes_R\Hom_R(E(R/\frak m),\cE))_n.$$  

We therefore have $$E(R/\frak m)\otimes_R^{\bf L}\RHom_R(D,G)\simeq \Sigma^t E(R/\frak m)\otimes_R\Hom_R(E(R/\frak m),\cE)$$$$\simeq \Sigma^t E(R/\frak m)\otimes_R\Hom_R(E(R/\frak m),G).$$

 Since $G$ is Gorenstein  injective, there exists an exact sequence of $R$-modules  $\dots\To E_1\To E_0\To G\To 0$ such that each $E_i$ is injective and $\Ker(E_i\To E_{i-1})$ is Gorenstein injective for each $i$. Setting $\cJ=\dots\To E_1\To E_0\To 0$, the complex $\Hom_R(E(R/\frak m),\cJ)$ is a flat resolution of $\Hom_R(E(R/\frak m),G)$. Hence we have the following isomorphism in $D(R)$
 $$E(R/\frak m)\otimes_R^{\bf L}\RHom_R(D,G)\simeq \Sigma^t(E(R/\frak m)\otimes_R \Hom_R(E(R/\frak m),\cJ));$$
 and consequently we have  $$\inf(E(R/\frak m)\otimes_R^{\bf L}\RHom_R(D,G))=\inf (\Sigma^t(E(R/\frak m)\otimes_R \Hom_R(E(R/\frak m),\cJ)))\geq t.$$    
\end{proof}

The following proposition generalizes the local duality theorem for arbitrary $R$-
modules.
\medskip
\begin{Proposition}\label{FS}
Let $\frak p$ be a prime ideal of $R$ such that $\height \frak p=h$ and let $M$ be an $R$-module. Then for every integer $k$, there is an isomorphism which is functorial in $M$
$$\Tor_k^R(M,\Hom_R(D,E(R/\frak p)))\cong H_{\frak pR_{\frak p}}^{h-k-\sup D_{\frak p}}(M_{\frak p}).$$
\end{Proposition}
\begin{proof}
As the functors $H_k(-\otimes_R^{\bf L}\Hom_R(D,E(R/\frak p)))$ and $H_{\frak pR_{\frak p}}^k(-)$ preserve the direct limits for each $k$, we may assume that $M$ is a finitely generated $R$-module. We also observe that $\Tor_k^R(M,\Hom_R(D,E(R/\frak p)))\cong \Tor_k^{R_{\frak p}}(M_{\frak p},\Hom_{R_{\frak p}}(D_{\frak p},E(R/\frak p)))$ for each $k$. Then without loss of generality, we may assume that $R$ is a local ring with the maximal ideal  $\frak p=\frak m$ and so $\dim R=h$. We notice that $\Sigma^{h-\sup D}D$ is a normalized dualizing complex (see [CF, Proposition 7.1.11]). Now using the local duality theorem [H, Corollary 6.3] for local cohomology modules, we have the following isomorphisms
$$\Tor_k^R(M,\Hom_R(D,E(R/\frak m)))\cong H_k(M\otimes_R^{\bf L}\Hom_R(D,E(R/\frak m)))$$$$\cong H_k(\Hom_R(\RHom(M,D),E(R/\frak m)))$$$$\cong H_{k-h+\sup D}(\Hom_R(\RHom(M,\Sigma^{h-\sup D}D),E(R/\frak m)))\cong   H_{\frak m}^{h-k-\sup D}(M).$$ The functionality of the isomorphism follows from the local duality theorem.
   \end{proof}

  For every integer $k$, we define $X_k=\{\frak p\in \Spec R|\hspace{0.1cm}{\rm ht}\frak p-\sup D_{\frak p}=k\}$.
\medskip
   
   \begin{Lemma}\label{lemGor}
   Let $G$ be a Gorenstein injective $R$-module and let $E$ be an injective $R$-module. Then for every integer $k$, there is an isomorphism which is functorial in $G$
   $$\Tor_k^R(E,\RHom_R(D,G))\cong \Tor_k^R(G,\Hom_R(D,E)).$$ 
   \end{Lemma}
   \begin{proof}
    
  By Matlis structure theorem, we have $E=\bigoplus_{\frak p\in\Spec R}E(R/\frak p)^{\mu_{\frak p}}$. As $D\in D_b^f(R)$, we have $\Hom_R(D,E)\simeq\bigoplus_{\frak p\in\Spec R}\Hom_R(D,E(R/\frak p))^{\mu_{\frak p}}$. Moreover, since the functor $\Tor_k^R(-,-)$ preserves direct sums, using \cref{secgor}, \cref{zero} and \cref{FS}, we may assume that $E=E(R/\frak p)$ for some $\frak p\in\Spec R$ and $k=\dim R_{\frak p}-\sup D_{\frak p}$. 
   Given an injective $R$-module $I$, we have $I=\bigoplus_{\frak q\in\Spec R}E(R/\frak q)^{\mu_{\frak q}}$ and so by \cref{FS} and the above argument, we have 
    $$\Tor _k^R(E(R/\frak p),\Hom_R(D,I))\cong \bigoplus_{\frak q\in\Spec R}\Tor _k^R(E(R/\frak p),\Hom_R(D,E(R/\frak q)))^{\mu_{\frak q}}$$$$\cong \bigoplus_{\frak q\in X_k}\Gamma_{\frak qR_{\frak q}}(E(R/\frak p)_{\frak q})^{\mu_{\frak q}}\cong \bigoplus_{\frak p\subset \frak q\in X_k}\Gamma_{\frak qR_{\frak q}}(E(R/\frak p)_{\frak q})^{\mu_{\frak q}}\cong \bigoplus_{\frak p\subset \frak q\in X_k}\Gamma_{\frak q}(E(R/\frak p))^{\mu_{\frak q}}.$$  If $\frak p\subsetneq \frak q$, then there exists $r\in\frak q\setminus\frak p$ an so $E(R/\frak p)\stackrel{r.}\To E(R/\frak p)$ is isomorphism so that $\Gamma_{\frak q}(E(R/\frak p))\stackrel{r.}\To \Gamma_{\frak q}(E(R/\frak p))$ is isomorphism. Since any element of $\Gamma_{\frak q}(E(R/\frak p))$ annihilated by some power of $r$, we conclude that $\Gamma_{\frak q}(E(R/\frak p))=0$. Hence $$\Tor _k^R(E(R/\frak p),\Hom_R(D,I))\stackrel{\theta^I}\cong\Gamma_{\frak pR_{\frak p}}(E(R/\frak p)^{\mu_{\frak p}}=E(R/\frak p)^{\mu_{\frak p}}$$ is injective. It follows from the local duality theorem that $\theta^{(-)}$ is functorial on injective modules. Since $G$ is Gorenstein injective, there exists an exact sequence of $R$-modules  $\dots\To E_1\To E_0\To G \to 0$ such that each $E_i$ is injective and $K_i=\Coker(E_{i+1}\To E_{i})$ is Gorenstein injective for each $i$ (consider $K_0=G$). Assuming $E_0=\bigoplus_{\frak q\in\Spec R}E(R/\frak q)^{\mu_{\frak a}}$, $E_1=\bigoplus_{\frak q\in\Spec R}E(R/\frak q)^{\nu_{\frak a}}$ and setting $F(-)=\Tor_k^R(E(R/\frak p),\RHom_R(D,-)$, we have the following commutative diagram with exact rows 
        
$$\xymatrix{F(E_1)\ar[d]^{\theta^{E_1}}\ar[r]& F(E_0)\ar[d]^{\theta^{E_0}}\ar[r]& F(G)\ar@{-->}[d]^{\theta^G}\ar[r]& 0\\
\Gamma_{\frak pR_{\frak p}}(E(R/\frak p))^{\mu_{\frak p}}\ar[r]& \Gamma_{\frak pR_{\frak p}}(E(R/\frak p))^{\nu_{\frak p}}\ar[r]& \Gamma_{\frak pR_{\frak p}}(G_{\frak p})\ar[r]& 0.}$$
The exactness of the top row is obtained from the Gorenstein injectivity of $K_i$s and \cref{zero}. Moreover, the exactness of the bottom row follows using \cref{secgor} and the fact that $K_{i_{\frak p}}$ is a Gorenstein injective $R_{\frak p}$-module for each $i$ (see [CFH1, Proposition 5.5]). The homomorphism $\theta^G$ exists by the universal property of cokernels. Now the five lemma implies that $\theta^G$ is isomorphism and $\theta^{(-)}$ is functorial on Gorenstein injective modules as it is functorial on injective modules. Finally, according to \cref{FS}, we have the following isomrphisms which are functorial in $G$ $$\Tor_k^R(E(R/\frak p),\RHom_R(D,G))\overset{\theta^G}\cong\Gamma_{\frak pR_{\frak p}}(G_{\frak p})\cong \Tor_k^R(G,\Hom_R(D,E(R/\frak p))).$$  
 \end{proof}
 
  The previous lemma can be extended to any $R$-module of finite Gorenstein injective dimension.  
   
   \medskip
   
   \begin{Theorem}\label{isom}
   Let $M$ be an $R$-module of finite Gorenstein injective dimension and let $E$ be an injective $R$-module. Then for every integer $k$, there is an isomrphism which is functorial in $M$ 
   $$\Tor _k^R(E,\RHom_R(D,M))\stackrel{\theta^M_k}\cong \Tor_k^R(M,\Hom_R(D,E)).$$
   \end{Theorem}
  \begin{proof}
 For any $\frak p\in\Spec R$ and any integer $k$, we have the following isomorpgisms 
 $$\Tor _k^R(E(R/\frak p),\RHom_R(D,M))\cong \Tor_k^{R_{\frak p}}(E(R/\frak p),\RHom_{R_{\frak p}}(D_{\frak p},M_{\frak p}));$$  $$\Tor_k^R(M,\Hom_R(D,E(R/\frak p)))\cong\Tor_k^{R_{\frak p}}(M,\Hom_{R_{\frak p}}(D_{\frak p},E(R/{\frak p})).$$  Moreover, the functor $\Tor_k^R(-,-)$ preserves direct sums for any integer $k$. Hence we may assume that $R$ is a local ring of dimension $d$ with the maximal ideal $\frak m$ and  $E=E(R/\frak m)$. Assume that $\Gid_RM=n$. By induction on $n$, we prove that $\Tor _k^R(E,\RHom_R(D,M))\cong H_{\frak m}^{d-k-\sup D}(M)$ for each $k$ and so the result follows by \cref{FS}. The case $n=0$ is clear by \cref{lemGor} and \cref{FS}. For $n>0$, there exists an exact sequence of $R$-modules $0\To M\To G\To N\To 0$ such that $G$ is Gorenstein injective and $\Gid_RN= n-1$. This gives rise to an exact triangle $M\To G\To N\To\Sigma M$ in $D(R)$. Applying the functors $E\otimes_R^{\bf L}\RHom_R(D,-)$ and $-\otimes_R^{\bf L}\Hom_R(D,E)$ to it, we have the following  exact triangles in $D(R)$
  $$E\otimes_R^{\bf L}\RHom_R(D,M)\To E\otimes_R^{\bf L}\RHom_R(D,G)\To E\otimes_R^{\bf L}\RHom_R(D,N)\To\Sigma E\otimes_R^{\bf L}\RHom_R(D,M)$$
 $$M\otimes_R^{\bf L}\Hom_R(D,E)\To G\otimes_R^{\bf L}\Hom_R(D,E)\To N\otimes_R^{\bf L}\Hom_R(D,E)\To\Sigma M\otimes_R^{\bf L}\Hom_R(D,E)$$
 
For convenience, set $S_k(-)=\Tor _k^R(E,\RHom_R(D,-))$ and $T_k(-)=\Tor_k^R(-,\Hom_R(D,E))$ for every  integer $k$. Computing the homology modules in the above exact triangles,  we have the following  diagram with exact rows 
  
  $$\xymatrix{S_{k+1}(G)\ar[r]\ar[d]^{\theta_{k+1}^G}&S_{k+1}(N)\ar[r]\ar[d]^{\theta_{k+1}^N}& S_k(M)\ar@{-->}[d]^{\theta_k^M}\ar[r]& S_k(G)\ar[d]^{\theta_k^G}\ar[r]& S_k(N)\ar[d]^{\theta_k^N}\\
T_{k+1}(G)\ar[r]&T_{k+1}(N)\ar[r]&T_k(M)\ar[r]& T_k(G)\ar[r]& T_k(N).}$$
The vertical maps $\theta_{k}^G,\theta_{k}^N, \theta_{k+1}^G,\theta_{k+1}^N$ and the commutativity of the left and the right squares follow by the induction hypothesis and \cref{lemGor} for each $k$. We notice that  
 $T_k(G)=S_k(G)=0$ for all $k\neq \dim R-\sup D$  by \cref{zero} and $T_k(X)=0$ for all $k>\dim R-\sup D$ and all $R$-modules $X$ be \cref{lemGor}. Hence the universal properties  of kernels and cokernels, together with the induction hypothesis ensure the existence of $\theta_k^M$ and the commutativity of the inner squares for each $k$. Using the induction hypothesis $\theta_{k+1}^G,\theta_{k+1}^N,\theta_k^G$ and $\theta_k^N$ are isomorphism and so the five lemma implies that $\theta_k^M$ is isomorphism. The functoriality of the isomorphism follows from \cref{nex}. 
   \end{proof} 

\medskip 
\begin{Lemma}\label{nex}
 The map $\theta^{(-)}_k:\Tor _k^R(E,\RHom_R(D,-))\To \Tor_k^R(-,\Hom_R(D,E))$ defined in \cref{isom}, is natural on the class of $R$-modules of finite Gorenstein injective dimension for every integer $k$.
\end{Lemma}
\begin{proof}
For every integer $k$, using Matlis structure theorem and the fact that the functor $\Tor_k^R(-,-)$ preserves direct sums, we may assume that $E=E(R/\frak p)$ for some $\frak p\in\Spec R$. For convenience, we set $S_k(-)=\Tor _k^R(E,\RHom_R(D,-))$ and $T_k(-)=\Tor_k^R(-,\Hom_R(D,E))$. It follows from \cref{isom} and \cref{zero} that $S_k(G)=T_k(G)=0$ for all $k\neq {\rm ht}\frak p-\sup D_{\frak p}$ and all Gorenstein injective $R$-modules $G$. Moreover, $S_k(M)=T_k(M)=0$ for all $k> {\rm ht}\frak p-\sup D_{\frak p}$ and all $R$-modules $M$ by \cref{FS}. Let $f:M_1\To M_2$ be a homomorphism of $R$-modules of finite Gorenstein injective dimension. Then by virtue of [Ho, Theorem 2.15], the $R$-modules $M_1$ and $M_2$ can be fitted into the exact sequences of $R$-modules $0\To M_1\To G_1\To N_1\To 0$ and $0\To M_2\To G_2\To N_2\To 0$ in which $G_1,G_2$ are Gorenstein injective $R$-modules and $N_1,N_2$ have finite injective dimensions. Hence we have the following commutative diagram of $R$-modules with exact rows
$$\xymatrix{0\ar[r]& M_1\ar[d]^f\ar[r]& G_1\ar[d]\ar[r]& N_1\ar[d]\ar[r]& 0\\
0\ar[r]& M_2\ar[r]& G_2\ar[r]& N_2\ar[r]& 0.}$$
We prove the claim by descending induction on $k$. If $k={\rm ht}\frak p-\sup D_{\frak p}$, we have the following  diagram with the exact rows 

$$\xymatrix@=0.5cm{
0\ar[rr]&&S_k(M_1)\ar[dl]_{S_k(f)}\ar[rrrr]\ar[dd]^{\theta^{M_1}_k} & & & & S_k(G_1)\ar[dl]\ar[dd]\\
0\ar[r]&S_k(M_2)\ar[rrrr]\ar[dd] & & & & S_k(G_2) \ar[dd] &\\
0\ar[rr]^{\theta^{M_2}_k}&& T_k(M_1)\ar[rrrr]\ar[dl]^{T_k(f)} & & & & T_k(G_1)\ar[dl]\\
0\ar[r]&T_k(M_2) \ar[rrrr] &&  & & T_k(G_2)  & \\}$$
We note that the five faces in the above diagram commute by \cref{lemGor} and the commutative diagram in \cref{isom}. Therefore $T_k(f)\theta^{M_1}_k=\theta_k^{M_2}S_k(f)$. For every $k<{\rm ht}\frak p-\sup D_{\frak p}$, by \cref{zero}, we have the following  diagram with the exact rows 
$$\xymatrix@=0.5cm{&&
S_{k+1}(N_1)\ar[dl]\ar[rrrr]\ar[dd] & & & & S_k(M_1)\ar[dl]^{S_k(f)}\ar[dd]^{\theta_k^{M_1}}\ar[r]&0\\
&S_{k+1}(N_2)\ar[rrrr]\ar[dd] & & & & S_k(M_2) \ar[dd]_{\theta_k^{M_2}} \ar[rr]&&0\\
&& T_{k+1}(N_1)\ar[rrrr]\ar[dl] & & & & T_k(M_1)\ar[dl]^{T_k(f)}\ar[r]&0\\
&T_{k+1}(N_2) \ar[rrrr] &&  & & T_k(M_2)  \ar[r]&0\\}$$ 
The induction hypothesis and the diagram in \cref{isom} imply that the five faces in the above diagram are commutative. Consequently $T_k(f)\theta_k^{M_1}=\theta_k^{M_2}S_k(f)$. 
\end{proof}

As an additional result, the isomorphisms in \cref{isom} can be generalized to any complex in $\cB(R)$.
     \begin{Proposition}\label{corcom}
   Let $M$ be a complex in $\cB(R)$ and let $E$ be an injective $R$-module. Then $$\Tor_k^R(E, \RHom_R(D,M))\cong\Tor_k^R(M,\Hom_R(D,E))$$
   for every integer $k$.
         \end{Proposition}
      \begin{proof}
   Since $M\in\cB(R)$, it has finite Gorenstein injective dimension by [CFH2, Theorem 10.4.7]. Assume that $\Gid_RM=n$. Without loss of generality, we may assume that $\sup M=0$ and $\cG:0\To G_0\To G_{-1}\To\dots\To G_{-n}\To 0$ is an Gorenstein injective resolution of $M$; that is $M\simeq \cG$ in $D(R)$. Replacing $M$ by $\cG$, we assume that $M=\cG$ is a bounded complex of Gorenstein injective modules. For any $\frak p\in\Spec R$, we have the following isomorphisms in $D(R)$ 
  $$E(R/\frak p)\otimes_R^{\bf L}\RHom_R(D,M)\simeq E(R/\frak p)\otimes_{R_{\frak p}}^{\bf L}\RHom_{R_{\frak p}}(D_{\frak p},M_{\frak p});$$  $$M\otimes_R^{\bf L}\Hom_R(D,E(R/\frak p))\simeq M_{\frak p}\otimes_{R_{\frak p}}^{\bf L}\Hom_{R_{\frak p}}(D_{\frak p},E(R/{\frak p})).$$ Hence, since the tensor preserves direct sums, we may assume that $R$ is a local ring with the maximal ideal $\frak m$; $E=E(R/\frak m)$ and $t=\dim R-\sup D$. For convenience, set  $S_k(-)=\Tor _k^R(E,\RHom_R(D,-))$ and $T_k(-)=\Tor_k^R(-,\Hom_R(D,E))$ for every integer $k$. We prove the claim by induction on $n$, The case $n=0$ was proved  in \cref{lemGor}. For $n=1$, by [H, Lemma 7.2], we have an exact triangle $G_0\To G_{-1}\To\cG\To\Sigma G_0$ in $D(R)$, where $\cG={\rm con}(I_0\To I_1)$. Applying $S_k(-)$ and $T_k(-)$ to the above exact triangle, it follows from \cref{zero} that $S_k(\cG)=T_k(\cG)=0$ for all $k\neq t, t+1$. Moreover, using \cref{zero}, we have the following commutative diagram with exact rows
  
  $$\xymatrix{0\ar[r]& S_{t+1}(\cG)\ar[r]\ar@{-->}[d]^{\theta_{t+1}^{\cG}}&S_{t}(G_0)\ar[r]\ar[d]^{\theta_{t}^{G_{0}}}& S_t(G_{-1})\ar[d]^{\theta_t^{G_{-1}}}\ar[r]& S_{t}(\cG)\ar@{-->}[d]^{\theta_{t}^{\cG}}\ar[r]& 0\\0\ar[r]&
T_{t+1}(\cG)\ar[r]&T_{t+1}(G_0)\ar[r]&T_t(G_{-1})\ar[r]& T_{t}(\cG)\ar[r]& 0}$$  
 where the middle square is commutative by \cref{lemGor}. The existence of $\theta_t^{\cG}$ and $\theta_{t+1}^{\cG}$, along with the fact that they are isomorphisms,  follows by universal properties of kernels and cokernels and the five lemma. Now, assume that $n\geq 2$. By virtue of [H, Lemma 7.2], there exists an exact triangle in $D(R)$ $$\cG_{>-n}\to \Sigma^{-n+1}G_{-n}\To \cG\To\Sigma \cG_{>-n}\hspace{0.4cm}(\dag_n),$$ where $\cG_{>-n}=0\to G_0\To\dots\To G_{-n+1}\To 0$ (we notice that $\cG={\rm con}(\cG_{>-n}\to \Sigma^{-n+1}G_n)$).
 
 We firs prove that $S_k(\cG)=T_k(\cG)=0$ for all $k$ such that: 
 
 \[
    k\neq
    \begin{cases}
    t-n+1,t-n+3,\dots,t-1,t+1,t+2,\dots,t+n; \hspace{0.2cm} {\rm if}\hspace{0.1cm} n \hspace{0.1cm}{\rm is\hspace{0.1cm} even}\\
  t-n+1,t-n+3,\dots,t,t+2,t+3,\dots,t+n; \hspace{0.2cm} {\rm if}\hspace{0.1cm} n \hspace{0.1cm}{\rm is\hspace{0.1cm} odd}.
 \end{cases}
 \]
Furthermore, we show that there exist the isomorphisms $S_k(\cG)\cong S_{k-1}(\cG_{>-n})$ and $T_k(\cG)\cong T_{k-1}(\cG_{>-n})$ for all $k=t-n+3,\dots, t+n$ in the above cases. Also $S_{t-n+1}(\cG)\cong S_{t}(G_{-n})$ and $T_{t-n+1}(\cG)\cong T_{t}(G_{-n})$. We prove the assertion for the functors $S_k(-)$  and the proof for functors $T_k(-)$ are the same using \cref{lemGor}.

 For the case $n=2$, applying $S_k(-)$ to $(\dag_2)$, we deduce that $S_k(\cG)=0$ for all $k< t-1$ and $k>t+2$. Furthermore, using \cref{zero}, we have the following exact sequences
  $$0\To S_{t+2}(\cG)\To S_{t+1}(\cG_{>-2})\To 0\To S_{t+1}(\cG)\To S_t(\cG_{>-2})\To 0.$$  
On the other hand, since Gid$_R\cG_{>-2}= 1$, using the case $n=1$, we have $S_k(\cG_{>-2})=0$ for all $k\neq t, t+1$ so that $S_{t-1}(\cG)\cong S_t(G_{-2}).$ Moreover, there exists an exact sequence $$0=S_{t+1}(G_{-2})\To S_t(\cG)\To S_{t-1}(\cG_{>-2})=0$$ which implies that $S_t(\cG)=0$.  For $n>2$, we assume that $n$ is odd and the same proof can be applied to the case where $n$ is even. For every integer $k$, application of the functor $S_k(-)$ to the exact triangle $(\dag_n)$ gives rise the following exact sequence of $R$-modules, 
$$S_{k}(\cG_{>-n})\To S_{k+n-1}(G_{-n})\To S_k(\cG)\To S_{k-1}(\cG_{>-n})\To S_{k+n-2}(G_{-n})\hspace{0.5cm}(\ddag).$$ 

We consider two cases. {\bf Case 1.} Vanishing of $S_k$: For $k=t-n+2, t-n+4,\dots,t+1$, since  $k\geq t-n+2$, we have $k+n-1\geq t+1$; and hence $S_{k+n-1}(G_{-n})=0$ by \cref{zero}. On the other, since $\Gid_R\cG_{>-n}= n-1$, using the induction hypothesis, we have $S_{k-1}(\cG_{>-n})=0$. Hence in view of the exact sequence $(\ddag)$, we have $S_k(\cG)=0$ for all $k=t-n+2, t-n+4,\dots,t+1$. For $k<t-n+1$ or $k>t+n$, using 
the induction hypothesis $S_{k-1}(\cG_{>-n})=0$ . Moreover, $S_{k+n-1}(G_{-n})=0$. Hence the exact sequence $(\ddag)$ forces that $S_k(\cG)=0$ for all  $k<t-n+1$ and $k>t+n$. {\bf Case 2.} Isomorphism of $S_k$: For every $k=t-n+3, t-n+5,\dots , t, t+2, t+3,\dots, t+n$, it follows from  \cref{zero} and the exact sequence $(\ddag)$ that $S_k(\cG)\cong S_{k-1}(\cG_{>-n})$.  Furthermore, using the induction hypothesis and the exact sequence $(\ddag)$, we have $S_{t-n+1}(\cG)\cong S_{t}(G_{-n})$. Therefore the assertion was proved. Now, using the induction hypothesis, we have 
$$ S_{t-n+1}(\cG)\cong S_{t}(G_{-n})\cong  T_{t}(G_{-n})\cong  T_{t-n+1}(\cG).$$ Moreover, for every $k=t-n+3, t-n+5,\dots ,t,t+2,t+3,\dots, t+n$, using the induction hypothesis, we have $$S_k(\cG)\cong S_{k-1}(\cG_{>-n})\cong T_{k-1}(\cG_{>-n})\cong T_k(\cG).$$ For all $k\neq t-n+1, t-n+3,\dots ,t,t+2,t+3,\dots, t+n$, we have $S_k(\cG)=T_k(\cG)=0$. For the case where $n$ is even, the proof is similarly.  
      \end{proof}
   
   As an immediate consequence, we establish the following result concerning the Gorenstein dimension of bounded complexes in $D(R)$.
   
   \begin{Corollary}\label{gorflat}
   Let $E$ be an injective $R$-module and let $M$ be a complex in $D_b(R)$. Then $$\Gid_RE\otimes_R^{\bf L}M\leq \Gfd_RM.$$ 
   \end{Corollary}
   \begin{proof}
   It suffices to consider that $M$ is a Gorenstein flat $R$-module and we show that $E\otimes_RM$ is a Gorestein injective $R$-module. It follows from [CFH1, Theorem 4.1] that $M\in\cA(R)$ and so we have
   $D\otimes_R^{\bf L}M\in\cB(R)$ so that $\Gid _R(D\otimes_R^{\bf L}M)$ is finite by [CFH1, Proposition 4.4]. Since $M\in\cA(R)$ and it is  Gorenstein flat, we have the following isomorphisms in $D(R)$
        $$E\otimes_RM\simeq E\otimes_R^{\bf L}M\simeq E\otimes_R^{\bf L}\RHom_R(D,D\otimes_R^{\bf L}M)).$$
 Therefore it follows from \cref{corcom} and the above quasi-isomorphisms that $$\Tor _0^R(D\otimes_R^{\bf L}M,\Hom_R(D,E))\cong E\otimes_RM$$ 
 and $\Tor _k^R(D\otimes_R^{\bf L}M,\Hom_R(D,E))=0$ for all $k>0$.  Since $D$ is dualizing complex, $\Hom_R(D,E)$ has finite flat dimension. Hence $(D\otimes_R^{\bf L}M)\otimes_R^{\bf L}\Hom_R(D,E)$ is a bounded complex of finite Gorenstein injective dimension. Assume that $\cG:0\To G_s\To G_{s-1}\To \dots \To G_t\To 0$ is a complex of Gorenstein injective modules such that  
  $\cG\simeq (D\otimes_R^{\bf L}M)\otimes_R^{\bf L}\Hom_R(D,E).$ By the previous argument, we have two exact sequences $0\To G_s\To\dots\To G_{-1}\To B_0\To 0$ and $\To Z_0\To G_0\To G_1\To \dots G_t\To 0$ which imply that $E\otimes_RM=Z_0/B_0$ has finite Gorenstein injective dimension. On the other hand , since $R$ has a dualizing complex, by [CF, Proposition 7.1.12], it has finite Krull dimension, say $d$. Thus it follows from [CFH1, Theorem 6.8] that $\Gid_R(E\otimes_RM)\leq d$. Since $M$ is Gorenstein flat, there exists a  flat resolution  $\dots F_1\To F_0\To 0$ of $M$  such that each $C_i=\Coker(F_{i+1}\To F_{i})$ is Gorenstein flat for each $i$. We notice that each $E\otimes_RF_i$ is injective and replacing $M$ by $C_d$, we conclude that $\Gid_R(E\otimes_RC_d)\leq d$. Consequently, $E\otimes_RM$ is a Gorenstein injective $R$-module.   
   \end{proof}
 
 \medskip
 The next result is crucial for constructing a finite filtration for Gorenstein injective modules.

   \begin{Proposition}\label{proo}
  Let $E$ be an injective $R$-module and let $G$ be a Gorenstein injective $R$-module. Then $\Tor_k^R(E,\RHom_R(D,G))$ is Gorenstein injective for every integer $k$. 
  \end{Proposition}
 \begin{proof}
As $E=\bigoplus_{\frak p\in\Spec R}E(R/\frak p)^{\mu_{\frak p}}$, for a fixed $k$ we have $$\Tor_k^R(E,\RHom_R(D,G))\cong \bigoplus_{\frak p\in\Spec R}\Tor_k^R(E(R/\frak p),\RHom_R(D,G))^{\mu_{\frak p}}$$$$\cong \bigoplus_{\frak p\in X_k}\Tor_k^R(E(R/\frak p),\RHom_R(D,G))^{\mu_{\frak p}}$$
 where the last isomorphisn is deduced from \cref{zero}. As $R$ has a dualzing module, by [CFH1, Theorem 6.9] direct sums of Gorenstein injective are Gorenstein injective. Hence by \cref{lemGor}, it suffices to prove that $\Tor_k^R(G,\Hom_R(D,E(R/\frak p))\cong \Tor _k^R(E(R/\frak p),\RHom_R(D,G))$ is Gorenstein injective where $k=\dim R_{\frak p}-\sup D_{\frak p}$. 
  Given an arbitrary injective module $I=\bigoplus_{\frak q\in\Spec R}E(R/\frak q)^{\mu_{\frak q}}$, by the same proof as mentioned in \cref{lemGor}, we have $\Tor _k^R(E(R/\frak p),\Hom_R(D,I))\cong E(R/\frak p)^{\mu_{\frak p}}$ and so is injective. Since $D$ has finite injective dimension, $\Hom_R(D,E(R/\frak p)$ has finite flat dimension. Thus there exists a complex of $R$-modules $\cF:=0\To F_m\To F_{m-1}\To\dots \To F_n\To 0$  such that each $F_i$ is flat and $\Hom_R(D,E(R/\frak p))\simeq \cF$. Application of  $G\otimes_R-$ to it gives rise to a complex of Gorenstein injective modules  $\cG=0\To G_m\To G_{m-1}\To \dots \To G_n\To 0$ such that  $G_i=G\otimes_RF_i$ for each $i$ and $G\otimes_R\Hom_R(D,E(R/\frak p)\simeq\cG$. Since $\Tor_i^R(G,\Hom_R(D,E(R/\frak p))\cong\Tor_i^R(E(R/\frak p),\RHom_R(D,G))=0$ for all $i\neq k$ by \cref{zero,lemGor}, the complexes $0\To G_m\To\dots\To G_{k+1}\To B_k\To 0$ and $0\To Z_k\To G_k\To \dots\To G_n\To 0$ are exact so that $B_k$ is Gorenstein injective and $Z_k$ has finite Gorenstein injective dimension. Therefore $\Tor_k^R(E(R/\frak p),\RHom_R(D,G))=Z_k/B_k$ has finite Gorenstein injective dimension. Now, since $G$ is is Gorenstein injective, there exist an exact sequence of $R$-modules $\dots\To E_1\To E_0\To G\To 0$ such that $E_i$ is injective and $K_i=\Coker(E_{i+1}\To E_i)$ is Gorenstein injective for each $i$. Setting $F(-)=\Tor _k^R(E(R/\frak p),\RHom_R(D,-))$ and replacing $G$ by $K_i$, the previous arguments deduce that $F(K_i)$ has finite Gorenstein injective dimension and $F(E_i)$ is injective for each $i$. On the other hand, using \cref{zero}, we have an exact sequence of modules $\dots\To F(E_1)\To F(E_0)\To F(G)\To 0$. Since $R$ has dualizing complex, dim$R$ is finite, say $d$. Then [CFH1, Theorem 6.8] implies that $\Gid_RF(K_{d-1})\leq d$ so that $F(G)$ is Gorenstein injective.   
\end{proof}
  
\medskip

By setting $\sup D=s$, for every $\frak p\in \Spec R$, we have the following (in)equalities  $$\id_RD\geq \id_{R_{\frak p}}D_{\frak p}=\id_{R_{\frak p}}\Sigma^{{\rm ht}\frak p-\sup D_{\frak p}}D_{\frak p}+{\rm ht}\frak p-\sup D_{\frak p}={\rm ht} {\frak p}-\sup D_{\frak p}\geq {\rm ht}\frak p-\sup D\geq -\sup D=-s$$
where the second equality holds as $\Sigma^{{\rm ht}\frak p-\sup D_{\frak p}}D_{\frak p}$  is a normalized dualizing complex. We find out from this description that there are only finitely many non-empty sets $X_k$; where $X_k=\{\frak p\in \Spec R|\hspace{0.1cm}{\rm ht}\frak p-\sup D_{\frak p}=k\}$. We observe that  $\Sigma^nD$ is a dualizing complex for every integer $n$ and so we may assume that $\sup D=0$. 
Now, we are prepared to present the main theorem of this section.  

\begin{Theorem}\label{fil1}
Let $\sup D=0$ and let $G$ be a Gorenstein injective module. Then $G$ has a finite filtration of Gorenstein injective submodules $$0=G_{d+1}\subset G_{d}\subset\dots \subset G_{1}\subset G_{0}=G$$ such that $G_{k}/G_{k+1}\cong \bigoplus_{\frak p\in X_{k}}\Tor_{k}^R(E(R/\frak p),\RHom_R(D,G))$ is Gorenstein injective for each $0\leq k\leq d$. Furthermore, such filtrations and direct sum decompositions are unique and functorial in $G$.
\end{Theorem}
\begin{proof}
It follows from [K, Proposition B.2] and [CF, Theorem 5.1.8]  that $D$ admits a minimal injective resolution $\cE:0\To E^0\To E^1\To\dots \To E^d\To 0$. On the other hand, using [CFH1, Theorem 3.3], we have $\inf\RHom_R(D,G)\geq 0$ and so assume that $\dots F_1\To F_0\To 0$ is a projective resolution of $\RHom_R(D,G)$. We notice that $D_{\frak p}$ is a dualizing complex for $R_{\frak p}$  with the minimal injective resolution $\cE_{\frak p}$ for each $\frak p\in\Spec R$. Then, by the same argument mentioned in the proof of \cref{zero}, for each $k$, we have $E^k=\bigoplus_{\frak p\in X_k}E(R/\frak p)$. Then we have the bicomplex
$M_{-p,q}=E_{-p}\otimes_RF_q$ which can be regarded as a first quadrant bicomplex shifting the indexes (we may set $E^p=E_{-p}$). The total complex induced by $(M_{-p,q})$ is Tot$(M)$ with the $n$-th term Tot$(M)_n=\bigoplus_{p=0}^d M_{-p,n+p}$ and since $G\in\cB(R)$, we have $H_0({\rm Tot}(M))\cong G$ and $H_n(({\rm Tot}(M))=0$ for all $n\neq 0$. The spectral sequence determined by the first filtration of Tot$(M)$ is $(E^r_{p,q}, d^r_{p,q})$ where $$E^1_{-p,q}=H_q(\dots\To E^p\otimes_RF_{q+1}\To E^p\otimes_RF_q\To E^p\otimes_RF_{q-1}\To \dots)$$$$=\Tor_q^R(E^p,\RHom_R(D,G))\cong\bigoplus_{\frak p\in X_p}\Tor_q^R(E(R/\frak p),\RHom_R(D,G)).$$
We observe that this spectral sequence converges to a graded module $H$ where $H_n=H_n({\rm Tot}(M))$ (see [Ro, Proposition 10.26]). In particular, there exists a finite filtration 
  $$0=\Phi^{d+1}H_0\subset \Phi^dH_0\subset\dots\subset \Phi^1H_0\subset \Phi^0H_0=H_0=G$$ 
where $\Phi^pH_0/\Phi^{p+1}H_0=E^{\infty}_{-p,p}$ for all $0\leq p\leq d$. It follows from \cref{zero} that $E_{-p,q}^1=0 $ for all $p\neq q$ and hence $E^2_{-p,p}=\Ker d^1_{-p,p}/\Im d^1_{-p+1,p}=E^1_{-p,p}$. Continuing in this manner,  we have $E_{-p,p}^{\infty}=E^1_{-p,p}= \bigoplus_{\frak p\in X_p}\Tor_p^R(E(R/\frak p),\RHom_R(D,G))$ is Gorenstein injective for all $0\leq p\leq d$ by \cref{proo}. By setting  $G_p=\Phi^pH_0$, the exact sequences $0\To G_{p+1}\To G_p\To G_p/G_{p+1}\To 0$ imply that $G_p$ is Gorenstein injective for  each $p$ as $G_d$ is Gorenstein injective. The functoriality and uniqueness of filtrations and direct decompositions are proved as in [EH, Theorem 3.1].  
\end{proof}
   
\section{A filtration via the section functors}
In this section, we assume that $\dim R=d$ and  
for any non-negative integer $k$, assume that $Y_k=\{\frak p\in\Spec R|\hspace{0.1cm} \height \frak p=k\}$.

\begin{Lemma}\label{sumgor}
Let $\frak a$ and $\frak b$ be ideals of $R$. If $G$ is a Gorenstein injective $R$-module, then so is $\Gamma_{\frak a}(G)+\Gamma_{\frak b}(G)$.
\end{Lemma}
\begin{proof}
It is clear that $\Gamma_{\frak a}(G)\cap\Gamma_{\frak b}(G)=\Gamma_{\frak a+\frak b}(G)$; and hence $\Gamma_{\frak a}(G)\cap\Gamma_{\frak b}(G)$ is Gorenstein injective by \cref{secgor}. Moreover, since $R$ has a dualizing complex, $\Gamma_{\frak a}(G)\oplus\Gamma_{\frak b}(G)$ is Gorenstein injective by [CFH1, Thorem 6.9]. Thus the assertion is concluded by the following exact sequence of modules
$$0\To \Gamma_{\frak a}(G)\cap\Gamma_{\frak b}(G)\To \Gamma_{\frak a}(G)\oplus\Gamma_{\frak b}(G)\To\Gamma_{\frak a}(G)+\Gamma_{\frak b}(G)\To 0.$$
\end{proof}

For every Gorenstein injective $G$, we now exhibit a finite filtration of Gorenstein submodules
$$0=G_{d+1}\subset G_d\subset G_{d-1}\subset\dots \subset G_1\subset G_0=G$$ 
by constructing the submodules  $G_k$ of $G$. This provides  an advantage over the filtration established in \cref{fil1}, where the submodules $G_k$ remained unknown.

\medskip
\begin{Theorem}\label{fil2}
Let $G$ be a Gorenstein injective module. Then $G$ has a finite filtration  of Gorenstein injective submodules $$0=G_{d+1}\subset G_d\subset G_{d-1}\subset\dots \subset G_1\subset G_0=G$$ such that $G_k= \Sigma_{\frak p\in Y_k}\Gamma_{\frak p}(G)$ and $G_k/G_{k+1}=\Sigma_{\frak p\in Y_k}\Gamma_{\frak p}(G/G_{k+1})$ for each $0\leq k\leq d$. Furthermore, such filtrations and sums are unique and functorial in $G$. 
\end{Theorem}
\begin{proof}
 Fix $0\leq k \leq d$ and set $G_J=\Sigma_J\Gamma_{\frak p}(G)$ where $J$ is an arbitrary finite subset of $Y_k$. It follows from \cref{sumgor} that $G_J$  is a Gorenstein injective submodule of $G$ and so $\underset{\To}{\rm lim}G_J=\Sigma_{\frak p\in Y_k} \Gamma_{\frak a}(G)$ is Gorenstein injective by [CFH1, Theorem 6.9].
We now set $G_k=\Sigma_{\frak p\in Y_k} \Gamma_{\frak a}(G)$ and we show that $G_{k+1}\subset G_k$ for each $0\leq k\leq d$. We first claim that $G_{k+1}=\Sigma_{\frak p\in Y_k}\Gamma_{\frak p}(G_{k+1})$. Since $G_{k+1}=\Sigma_{\frak q\in Y_{k+1}}\Gamma_{\frak q}(G)$, for every $\frak q\in Y_{k+1}$, there exists $\frak p\in Y_k$ such that $\frak p\subset\frak q$. Hence $\Gamma_{\frak q}(G)= \Gamma_{\frak p}(\Gamma_{\frak q}(G))\subseteq\Gamma_{\frak p}(G_{k+1})$ and consequently $G_{k+1}\subset\Sigma_{\frak p\in Y_k}\Gamma_{\frak p}(G_{k+1})$. The other side is clear. Therefore, we have the following isomorphism 
 $$G_{k+1}=\Sigma_{\frak p\in Y_k}\Gamma_{\frak p}(G_{k+1})\subset\Sigma_{\frak p\in Y_k}\Gamma_{\frak p}(G)=G_k.$$  According to \cref{secgor}, for each $\frak p\in Y_k$, application of the functor $\Gamma_{\frak p}$ to the exact sequence $0\To G_{k+1}\To G\To G/G_{k+1}\To 0$ gives an isomorphism  $\Gamma_{\frak p}(G/G_{k+1})\cong \Gamma_{\frak p}(G)/\Gamma_{\frak p}(G_{k+1})$. Therefore, $\bigoplus_{\frak p\in Y_k}\Gamma_{\frak p}(G)/\bigoplus_{\frak p\in Y_k}\Gamma_{\frak p}(G_{k+1})\cong \bigoplus_{\frak p\in Y_k}(\Gamma_{\frak p}(G)/\Gamma_{\frak p}(G_{k+1}))\cong \bigoplus_{\frak p\in Y_k}\Gamma_{\frak p}(G/G_{k+1})$. For each $0\leq k\leq d$, there exists an exact sequence of $R$-modules with exact rows and injective vertical maps 
  $$\xymatrix{0\ar[r]&U\ar[r]\ar[d]& \bigoplus_{\frak p\in Y_k}\Gamma_{\frak p}(G_{k+1})\ar[d]\ar[r]&G_{k+1}\ar[d]\ar[r]& 0\\
0\ar[r]&V\ar[r]&\bigoplus_{\frak p\in Y_k}\Gamma_{\frak p}(G) \ar[r]& G_k\ar[r]& 0.}$$
 It follows from the above argument and the snake lemma that $G_k/G_{k+1}=\Sigma_{\frak p\in Y_k}\Gamma_{\frak p}(G/G_{k+1})$.  Finally, we assert that $G_0=G$. Otherwise, there exists a prime ideal $\frak p\in\Ass_R(G/G_0)$ and so $\Gamma_{\frak p}(G/G_0)\neq 0$. On the other hand, there exists a minimal prime ideal $\frak q$ of $R$ such that $\frak q\subset \frak p$. By the definition $\Gamma_{\frak p}(G)\subset \Gamma_{\frak q}(G)\subset G_0$; and hence the exact sequence $0\To G_0/\Gamma_{\frak p}(G)\To G/\Gamma_{\frak p}(G)\To G/G_0\To 0$ of Gorenstein injective modules and \cref{secgor} imply that $\Gamma_{\frak p}(G/G_0)=0$ which is a contradiction. 
  To prove the second claim, if $H$ is another Gorenstein injective $R$-module with  such a filtration $0=H_{d+1}\subset H_d\subset \dots\subset H_1\To H_0=H$ and $f:G\To H$ is a $R$-homomorphism, then for each $1\leq k\leq d$ and $\frak p\in \Spec R$ with $\height \frak p=k$, we have an $R$-homomorphism $\Gamma_{\frak p}(f):\Gamma_{\frak p}(G)\To \Gamma_{\frak p}(H)$ (the restriction of $f$) which gives rise to a restriction homomorphism  $f_k=\Sigma_{\frak p\in Y_k}\Gamma_{\frak p}(f):G_k\To H_k$   and the following commutative diagram 
 $$\xymatrix{0\ar[r]&G_{k+1}\ar[r]\ar[d]^{f_{k+1}}& G\ar[d]^f\ar[r]& G/G_{k+1}\ar[d]^{\bar{f}_{k+1}}\ar[r]& 0\\
0\ar[r]&H_{k+1}\ar[r]&H\ar[r]& H/H_{k+1}\ar[r]& 0.}$$  
Using the universal properties of coproducts  and sums, we have the following homomorphism $$\Sigma_{\frak p\in Y_k}\Gamma_{\frak p}(\bar{f}_{k+1}):G_k/G_{k+1}\To H_k/H_{k+1}.$$ 
\end{proof}

\medskip
Over a Gorenstein ring, Enochs and Huang [EH, Remark 3.3] remarked that the referee pointed out that the $G_k$ appearing in \cref{fil1} can be described via the formula $G_k/G_{k+1}=\bigoplus_{\frak p\in Y_k}\Gamma_{\frak p}(G/G_{k+1})$ for every integer $k$. However, they did not know how this holds. In the following theorem, we  establish this formula in a more general, where $R$ is assumed to admit a dualizing complex. 
\begin{Theorem}
Let $G$ be a Gorenstein injective $R$-module with the the same filtration  as in \cref{fil2}. Then $G_k/G_{k+1}=\bigoplus_{\frak p\in Y_k}\Gamma_{\frak p}(G/G_{k+1})$ for every integer $k$. Furthermore, the filtrations in \cref{fil1,fil2} are the same.
\end{Theorem}
\begin{proof}
With the same notations as  in \cref{fil2}, we prove that the sum occurring in $G_k/G_{k+1}$ is the direct sum for each $k$. Fixing an integer $k$, according to \cref{fil2}, we have $G_k/G_{k+1}=\Sigma_{\frak p\in Y_k}\Gamma_{\frak p}(G/G_{k+1})=\underset{\To}{\rm lim}G_J$, where $G_J=\Sigma_J\Gamma_{\frak p}(G/G_{k+1})$ and  $J$ is an arbitrary finite subset of $Y_k$. Then for every finite subset $J=\{\frak p_1,\dots,\frak p_n\}$ of $Y_k$, it suffices to show that $\Sigma_J\Gamma_{\frak p}(G/G_{k+1})=\bigoplus_J\Gamma_{\frak p}(G/G_{k+1})$. We proceed by induction on $n$. The case $n=1$ is clear. Assume that $n>1$.
 With re-indexing and using the induction hypothesis, it is enough to prove that 
$\Gamma_{\frak p_1}(G/G_{k+1})\bigcap\bigoplus_{J\setminus\{\frak p_1\}}\Gamma_{\frak p}(G/G_{k+1})=0$. To do this, if $x\in\Gamma_{\frak p_1}(G/G_{k+1})\bigcap\bigoplus_{J\setminus\{\frak p_1\}}\Gamma_{\frak p}(G/G_{k+1})$, then $x=x_2+\dots++ x_n$ where $x_i\in\Gamma_{\frak p_i}(G/G_{k+1})$. Thus, $x_i\in\Gamma_{\frak p+\frak p_i}(G/G_{k+1})$ for every $2\leq i\leq n$.  Fixing $i$, it is clear that ${\rm ht}(\frak p+\frak p_i)\geq k+1$ and so we have $\sqrt{\frak p+\frak p_i}=\bigcap_{i=1}^t\frak q_i$ where ${\rm ht}\frak q_j\geq k+1$ for each $j$. Thus $\Gamma_{\frak q_j}(G)\subset G_{k+1}$ for each $j$, and so there exists an exact sequence of Gorenstein injective modules 
$$0\To G_{k+1}/\Gamma_{\frak q_j}(G)\To G/\Gamma_{\frak q_j}(G)\To G/G_{k+1}\To 0.$$
Applying the functor $\Gamma_{\frak q_j}(-)$ to it and using \cref{secgor}, we deduce that $\Gamma_{\frak q_j}( G/G_{k+1})=0$ for each $j$. We prove by induction on $t$ that $\Gamma_{\frak p+\frak p_i}(G/G_{k+1})=\Gamma_{\cap_{i=1}^t\frak q_i}(G/G_{k+1})=0$. The case $t=1$ was proved. For $t>1$, by the Mayer-Vietoris sequence [BS, Theorem 3.2.3] and \cref{secgor}, there exists a short exact sequence of $R$-modules 
$$0\To\Gamma_{\frak q_1+\cap_{i=2}^t\frak q_i}(G/G_{k+1})\To \Gamma_{\frak q_1}(G/G_{k+1})\oplus \Gamma_{\cap_{i=2}^t\frak q_i}(G/G_{k+1})\To\Gamma_{\cap_{i=1}^t\frak q_i}(G/G_{k+1})\To 0.$$
The induction hypothesis implies that $\Gamma_{\frak q_1}(G/G_{k+1})=\Gamma_{\cap_{i=2}^t\frak q_i}(G/G_{k+1})=0$ and so the above exact sequence forces that $\Gamma_{\cap_{i=1}^t\frak q_i}(G/G_{k+1})=0.$  Hence $x_i=0$ for all $2\leq i\leq n$ and so $x=0$. Therefore $G_k/G_{k+1}=\bigoplus_{\frak p\in Y_k}\Gamma_{\frak p}(G/G_{k+1})$ for each $k$. We observe  by construction that the filtrations in \cref{fil1} and \cref{fil2} have the same length. Consequently, the uniqueness of such filtrations for $G$ ensures that they coincide.
\end{proof}

\begin{Remark}
A natural question posed by Feickert and Sather-Wagstaff [FS, Question 5.4] asks that for a Gorenstein injective $R$-module $G$, does there exist a decomposition $G\cong\oplus_{\lambda}G_{\lambda}$, where each $G_{\lambda}$ is Gorenstein injective? However, [FS, Example 5.3] shows that the filtrations arising from \cref{fil1} and \cref{fil2} need not yield such a direct sum decomposition $G\cong \bigoplus_k G_{k}/G_{k+1}$, even when $R$ is a Gorenstein local ring of dimension one. We also remark that an indecomposable Gorenstein injective $R$-module may admit a non trivial filtration (i.e. a filtration of length $\geq 1$).
\end{Remark}

\medskip

\begin{Example}[FS, Example 5.3]
Let $k$ be a field and set $R:=k[[X,Y]]/(X^2)$ with  maximal ideal $\frak m:=(X,Y)R$ and $E:=E(k)$. Set $\frak q:=(X)R$ and $\overline{R}:=R/\frak q$. We observe that $\Spec R=\{\frak q,\frak m\}$. Since $R$ is a Gorenstein ring of dimension one,  $G:=E_R(\overline{R})/\overline{R}$ is a Gorenstein injective $R$-module by [EJ, Thorem 10.1.13]. Furthermore, $G$ is indecomposable by [FS, Example 5.3]. It follows from \cref{fil2} that $G$ has a filtration $0=G_2\subset G_1\subset G_0=G$ of Gorenstein injective submodules in which $G_1=\Gamma_{\frak m}(G)$ and $G=\Gamma_{\frak q}(G)$ (we notice that $G$ is $\frak q$-torsion). We observe that $G\ncong G_1\oplus G/G_1$ as $G$ is indecomposable. 
\end{Example}


\end{document}